\documentclass[preprint,12pt]{amsart}

\usepackage{amsmath,amsfonts,amsthm,amssymb,calrsfs,mathrsfs,enumerate}

\usepackage[a4paper,twoside,left=3cm,right=3cm,bottom=3cm,top=3cm]{geometry}

\usepackage[all,2cell, dvips,pdftex]{xy} \SilentMatrices
\usepackage{epsfig,color}

\usepackage{bm}
\usepackage{pst-plot}

\theoremstyle{plain}
\newtheorem{theorem}{Theorem}[section]

\newtheorem*{lemma*}{Lemma}
\newtheorem{lemma}[theorem]{Lemma}

\newtheorem*{theorem*}{Theorem}

\newtheorem{proposition}[theorem]{Proposition}
\newtheorem*{proposition*}{Proposition}

\newtheorem*{corollary*}{Corollary}

\theoremstyle{definition}

\newtheorem*{remark*}{Remark}

\newtheorem*{definition*}{Definition}

\newtheorem*{example*}{Example}

\def\gq{/\!\! /}

\def\RR{{\mathbf R}}

\def\Spec{\operatorname{Spec}}

\def\Cox{\operatorname{Cox}}

\def\ZZ{{\mathbf Z}}

\def\PP{{\mathbf P}}

\def\G{{\mathscr G}}

\def\O{{\mathscr O}}

\def\K{{\mathscr K}}

\def\+{\oplus}                   
\def\*{\otimes}                  
\def\hpil{\longrightarrow}       

\def\Eff{\operatorname{Eff}}

\def\Pic{\operatorname{Pic}}
\def\Hom{\operatorname{Hom}}

\def\f{\Gamma}

\begin{document}



\title{Cox rings of K3 surfaces with Picard number two }
\author{John Christian Ottem}


\address{Matematisk institutt, Universitetet i Oslo, PO Box 1053, Blindern, NO-0316 Oslo, Norway}


\begin{abstract}
We study generators and relations of Cox rings of K3 surfaces of Picard number two. In particular we consider the Cox rings of classical examples of K3 surfaces, such as quartic surfaces containing a line and elliptic K3 surfaces.
\end{abstract}




\maketitle

\thispagestyle{empty}

\section{Introduction}
\noindent Linear systems and projective models of K3 surfaces are classical objects in algebraic geometry and have been studied by many authors, see e.g. \cite{Kov94}, \cite{SD74} and \cite{TJ04}. In this paper we make use of these results to study explicitly the Cox ring of a K3 surface, which is essentially defined as the ring of all sections of all line bundles on the surface. It's well-known that the generic K3 surface $X$ has Picard number one, so in this case the Cox ring is just the section ring $\bigoplus_{m\ge 0}H^0(X,mD)$, where $D$ is an ample divisor generating $\Pic(X)$. Generators and relations of this ring were investigated by Saint-Donat in \cite{SD74}. In this paper we will consider the case when $X$ has Picard number two.

Cox rings of K3 surfaces were first studied by Artebani, Hausen and Laface in the recent paper \cite{AHL09}. In that paper it was shown that a K3 surface has a finitely generated Cox ring if and only if its effective cone is rational polyhedral. When the Picard number is two, it is known that this cone is rational polyhedral if an only if $\Pic(X)$ contains a class of self-intersection 0 or $-2$. On the other hand, for higher Picard number, the cone is rational polyhedral only when $\Pic(X)$ belongs to some finite list of hyperbolic lattices (see \cite[Theorem 2.12]{AHL09}). This makes the case when $X$ has Picard number two particularly interesting. 

Knowing when the Cox ring is finitely generated raises the problem of finding an explicit presentation for it. This problem was considered in \cite{AHL09} for many classes of K3 surfaces, including double covers of Del Pezzo surfaces. The aim of this paper is to extend some of these results and also present some new examples.

The paper is organized as follows: In Section 2 we give a new proof of finite generation of the Cox ring when the effective cone is rational polyhedral. The results of this section can be used to study Cox rings of any Picard number. Our proof is also constructive in the sense that we provide a set of generators for $\Cox(X)$ which can be used to find a presentation of the Cox ring. The main idea here is to use Koszul cohomology to study the multiplication maps on global sections of line bundles on $X$ which gives the ring structure on $\Cox(X)$. In Section 3 we consider three examples where Theorem \ref{generators} is used to compute the Cox rings explicitly. In Section \ref{2ell}, we study in detail K3 surfaces with intersection matrix of the form $\left(\begin{smallmatrix}0 & d\\ d & 0 \end{smallmatrix}\right)$. Using classical results on elliptic fibrations in rational normal scrolls we are able to study the minimal resolution of the Cox ring. Our main result here is that the Betti numbers of this resolution coincide with the Betti numbers of elliptic normal curves.

\section{Cox rings of K3 Surfaces}

\noindent Let $X$ be a smooth projective K3 surface with $\Pic(X)$ freely generated by the classes of  effective divisors $\f_1,\ldots,\f_{\rho}$. The Cox ring of $X$ is defined by
$$\Cox(X)=\bigoplus_{{\bf m}\in \mathbb{Z}^\rho}H^0\left(X,m_1\f_1+\ldots+m_{\rho}\f_{\rho}\right)$$where the ring structure is given by viewing sections as rational functions. Note that $\Cox(X)$ comes with a multigrading given by the semigroup of effective divisors. By Theorem 2.7 in \cite{AHL09} the Cox ring is a finitely generated $k$-algebra if and only if this semigroup is finitely generated, or equivalently, the cone of effective divisors $\Eff(X)\subset \Pic(X)\otimes \RR$ is rational polyhedral. In this case, a result of Kov\'acs \cite{Kov94} says that this cone is spanned by classes of self-intersection $-2$ or $0$.

To find a concrete presentation of $\Cox(X)$, we look for a minimal set of generators $x_1,\ldots,x_r$ from respective vector spaces $H^0(X,D_1),\ldots, H^0(X,D_r)$ and describe the ideal of relations $I$ between them. Letting $R=k[x_1,\ldots,x_r]$ with the natural $\Pic(X)$-grading given by $\deg x_i=D_i$, there is an exact sequence of $\Pic(X)$-graded $k$-algebras
$$0\to I\to R\to \Cox(X)\to 0.
$$As explained in \cite{HK00}, this presentation gives an embedding of $X$ into a toric variety. Indeed, consider the affine variety $V(I)\subset \mathbb{A}^r=\Spec(R)$. There is a natural action of the Neron-Severi torus $G=\Hom(\Pic(X),\mathbb{G}_m)\simeq \mathbb{G}_m^\rho$ on $R$ and $\Cox(X)$, making the inclusion $V(I) \subset \mathbb{A}^r$ $G$-equivariant. Taking GIT quotients, we get an inclusion  $i:X=V(I)\gq G \hookrightarrow\mathbb{A}^r\gq G$. The quotient $Y=\mathbb{A}^r \gq G$ is a normal toric variety with $R$ as its Cox ring. By \cite[Prop. 2.11]{HK00}, this embedding induces an isomorphism $i^*:\Pic(Y)\to \Pic(X)$, taking the effective cone of $Y$ to the effective cone of $X$. This toric embedding will be useful in Section \ref{2ell}. Also, in the case we are interested in, the torus $G$ is 2-dimensional, and so it is easy to see that the Krull dimension of the Cox ring is 4 (for a complete proof see \cite{BP04}).


%

We recall a few standard facts about linear systems on K3 surfaces:

\begin{proposition}\label{nefp}
Let $X$ be a smooth projective K3 surface, $D\neq 0$ be an effective divisor on $X$.
\begin{enumerate}[i)]
\item If $D$ is nef, then the linear system $|D|$ has a base-point if and only if there exist curves $E$, $\f$ and an integer $k\ge 2$ such that
\begin{equation}\label{fixedcomp}
D\equiv kE+\f, \quad E^2=0, \quad \f^2=-2, \quad E\cdot \f=1.
\end{equation}
\item If $D^2\ge 0$, then $H^1(X,D)\neq 0$ if and only if either $i)$ $D\equiv kE$ for some divisor $E$ with $E^2=0$ and $k\ge 2$ or $ii)$ $D\cdot \f\le -2$ for some divisor $\f$ with $\f^2=-2$.
\item If $D$ is nef and $D^2=0$ then $D$ is base-point free and $D\equiv kE$ for some smooth elliptic curve $E$. If $D^2>0$, then the generic element of $|D|$ is a smooth and irreducible curve.
\item A nef divisor class $D$ with $D^2>0$ is hyperelliptic (i.e., the generic member of $|D|$ is a hyperelliptic curve) if and only either $D^2=2$; or there is a smooth elliptic curve $E$ such that $D\cdot E=2$; or $D=2B$ for a smooth curve $B$ with $B^2=2$
\item If $D$ is not hyperelliptic, then the section ring $R(X,D)=\bigoplus_{n\ge 0}H^0(X,nD)$  is generated in degree 1. If $D$ is hyperelliptic and $g(D)=2$, then $R(X,D)$ is generated in degree 1 and 3. If $D$ is hyperelliptic and if $g(D)\ge 3$, then $R(X,D)$ is generated in degrees 1 and 2. 
\end{enumerate}
\end{proposition}

\begin{proof}
i) follows from \cite[2.7]{SD74} and ii) from the main result of \cite{KL07}. Then $iii)$ and $iv)$ follow from \cite[Proposition 2.6]{SD74}, while $v)$ follows from \cite{SD74} and \cite[Proposition 3.4]{AHL09}
\end{proof}

We prove the following result on the surjectivity of the multiplication maps on a K3 surface:

\begin{proposition}\label{koszulK3}
Let $X$ be a smooth projective K3 surface. Let $D$ and $E$ be nef divisors on $X$ such that $|E|$ is base-point free. Then the multiplication map
\begin{equation}\label{multmap}
H^0(X,D)\otimes H^0(X,E)\to H^0(X,D+E)
\end{equation}is surjective if $H^1(X,D-E)=H^1(X,D)=0$ and $H^2(X,D-2E)=0$.
\end{proposition}


\begin{proof}Proving the lemma is equivalent to showing that the Koszul cohomology group $\K_{0,1}(X,D,E)$ is zero (see \cite{MG84}). By definition, this is the homology of the complex$$
\bigwedge^1 H^0(X,E)\otimes H^0(X,D)\to \bigwedge^0 H^0(X,E)\otimes H^0(X,D+E)\to 0.
$$Now, the assumptions on the cohomology vanishing ensure us that we are in position to apply the duality theorem of \cite{MG84}, which states that under these circumstances, $
 \K_{0,1}(X,D,E)\cong \K_{r-2,2}(X,-D,E)^*
$ where $r=h^0(X,E)-1$ and $\K_{r-2,3}(X,-D,N)$ is the homology of the complex
\begin{eqnarray*}
\bigwedge^{r-1} H^0(X,E)\otimes H^0(X,-D+E)&\to& \bigwedge^{r-2} H^0(X,E)\otimes H^0(X,-D+2E)\\ &\to& \bigwedge^{r-3} H^0(X,E)\otimes H^0(X,-D+3E).
\end{eqnarray*}But by assumption, $H^0(X,-D+2E)=0$ and so the homology of the complex is zero.\end{proof}

We now define a set $\G$ that contains the degrees of the generators of $\Cox(X)$. 

 
 \begin{enumerate}[A)]
\item Let $\G$ be the set of all classes of curves with self-intersection $-2$.
 \item Add to $\G$ the nef divisors $D$ such that for every base-point-free divisor $E$ with $D-E$ effective, either $H^1(X,D-2E)\neq 0$;  or $H^1(X,D-E)\neq 0$; or $H^2(X,D-3E)\neq 0$.
 \item If $D$ from $B)$ was non-hyperelliptic, remove all higher multiples $mD$ $m\ge 2$ from $\G$.
 \item If $D$ from $B)$ was hyperelliptic, remove $2D$ from $\G$ if $D^2=2$ and $3D$ if $D^2>2$.
 \end{enumerate}
 
\noindent Some remarks are in order here. First, we will prove in Theorem \ref{generators} below that $\G$ is actually a finite set. Moreover, we point out that thanks to the classical results in Proposition \ref{nefp}, finding the set $\G$ is straightforward once one has a description of the effective cone. Indeed, $(-2)$-curves in $A)$ are extremal in the effective cone and can be found by inspection. Furthermore,  the nef divisors satisfying $B)$ can be obtained easily using Proposition \ref{nefp}ii) and the observation that $H^2(X,D-3E)\neq 0$ if and only if $3E-D$ is effective. The examples in Section 3 will make this clear. Finally, the conditions $C)$ and $D)$ are there to eliminate the redundant degrees of $R(X,D)$, in accordance with Proposition \ref{nefp}v). 

 \begin{theorem}\label{generators}
 Let $X$ be a smooth projective K3 surface with rational polyhedral effective cone and let $\G$ be defined by $A)-D)$  above. The $\G$ is finite and the Cox ring of $X$ is finitely generated by sections of degrees contained in $\G$.
 \end{theorem}
 
 \begin{proof}
We first show that the set $\G$ is finite. First, the classes of curves of self-intersection $-2$ are extremal in the effective cone, so if this cone is assumed to be rational polyhedral, we see that the classes satisfying $A)$ is finite.

We show that the set of nef divisors $D$ satisfying $B)$ and $C)$ above is finite. We may assume $D^2>0$ since if $D^2=0$, $D$ is linearly equivalent to $kE$ for some elliptic curve $E$ and in that case $D\not \in \G$ if $k\ge 2$, by $C)$ above (also there are only finitely many such $E$ since these are extremal in the nef cone). We will need the fact that on a K3 surface, $|3N|$ is base-point free for $N$ a nef divisor \cite{SD74}. Then if $N$ runs over the set of non-zero nef divisors, the union of the translates $10N+\mbox{Nef}(X)$ covers all but finitely many integral nef divisor classes. So it suffices to show that if $D$ is a nef divisor class contained in this union, then $D\not\in \G$. So suppose $D$ is a nef divisor of the form $10N+D'$ with $N,D'$ nef. Then if we let $E=3N$, we see that $D-kE=N+D'$ is nef and big for $k=1,2,3$ and so $D$ and $E=3N$ satisfy $B)$, and hence $D\not \in \G$. In all, this shows that $\G$ is finite.

We now show that the Cox ring of $X$ is finitely generated by sections of degrees contained in $\G$. Fix a very ample divisor $H$ on $X$, so that we may talk about the \emph{degree}, $H\cdot D>0$, of an effective divisor class $D$. Let $D$ be an effective divisor class. We show that any section $s\in H^0(X,D)$ can be written as a polynomial in the sections in the above degrees using Proposition \ref{koszulK3} and induction on the degree. 

 If $D$ is not base-point free, then by Proposition \ref{nefp} there is a $(-2)$-curve curve $\f$ (hence $\f\in \G$) in the  base locus of $|D|$. If $x$ is a section defining $\f$ then multiplication by $x$ gives an isomorphism $H^0(X,D-\f)\xrightarrow{\cdot x} H^0(X,D)$. Hence every section of $H^0(X,D)$ can be written as a product of $x\in H^0(\f)$ and an element of $H^0(X,D-\f)$. Replacing $D$ by $D-\f$, we may therefore reduce to the case where $D$ is base-point free.
 
Now $\G$ is constructed such that if $D$ is a base-point free divisor which is not in $\G$ then there is a base-point-free divisor $E$ such that the multiplication map \begin{equation}\label{mult}H^0(X,D-E)\otimes H^0(X,E)\to H^0(X,D).\end{equation}is surjective. By induction, elements of $H^0(X,E)$ and $H^0(X,D-E)$ are generated by sections of degrees contained in $\G$ and hence so the same applies to $H^0(X,D)$. \end{proof}

In particular, if the Picard number of $X$ is 2, the Cox ring is finitely generated if and only if $\Pic(X)$ contains a class of self-intersection $-2$ or $0$.

 \section{Examples}
 
 In this section we will demonstrate how Theorem \ref{generators} can be used to find explicit presentations of Cox rings, provided that the defining ideal is not too complicated.
 
 \subsection{A quartic surface containing a line}

We consider a quartic K3 surface $X$ with $\Pic(X)=\ZZ\f_1\oplus \ZZ\f_2$ with the intersection matrix given by $(\f_1\cdot \f_2)=\left(\begin{smallmatrix}-2 & \,\,3 \\\,\,3 & 0 \\ \end{smallmatrix}\right)$. Here we fix smooth curves $\f_1$ and $\f_2$ of genus 0 and 1 respectively. Such surfaces were studied in \cite{GM00}, where the authors refer to them as the \emph{Mori quartics}. Indeed, in this case it is straightforward to check that the divisor $H=\f_1+\f_2$ embeds $X$ as a quartic surface in $\PP^3$ and that $\f_1$ is sent to a projective line under this embedding.
It is also not hard to check that the effective cone is generated by the classes of $\f_1$ and $\f_2$ and that the nef cone is spanned by $\f_2$ and $3\f_1+2\f_2$ (see e.g., \cite[Proposition 3.1]{AHL09}).

The diophantine equation $(x\f_1+y\f_2)^2=-2x^2+6xy=-2$ has $(\pm 1,0)$ as the only solutions and so the only $(-2)$-curve on $X$ is $\f_1$. Moreover, using Proposition \ref{nefp}$iv)$, we find that there are no hyperelliptic classes on $X$. In particular, by Proposition \ref{nefp}$v)$ the section rings $R(X,D)$ are all generated in degree 1. Also, we find that for every ample divisor $D$ except the classes in $\G$ we have $H^1(X,D-2E)=H^1(X,D-E)=0$ and $H^2(X,D-3E)=0$ for some $E\in \{\f_2,\f_1+\f_2\}$. In the notation of Theorem \ref{generators}, this means that $$\G=\{\f_1,\f_2,\f_1+\f_2,3\f_1+2\f_2 \}$$ and using Riemann-Roch we find that we need sections $x,y_1,y_2,z_1,z_2,t$ such that 
$$
\begin{array}{ll}
H^0(X,\f_1)=\langle x \rangle & H^0(X,\f_1+\f_2)=\langle xy_1,xy_2,z_1,z_2\rangle\\
H^0(X,\f_2)=\langle y_1,y_2\rangle &H^0(X,3\f_1+2\f_2)=\langle x^3y_1^2,x^3y_1y_2,\ldots,z_2^2y_2,t\rangle
\end{array}
$$

For computing the defining ideal of the Cox ring we will need the following trick: If $R$ is a ring and $I\subset R[t]$ is an ideal containing an element of the form $ty+f$ with $f\in R$ and $y$ a non-zero divisors modulo $I$, then $I$ is prime if and only if the elimination ideal $I\cap R$ is prime. This can be seen by localizing $R[t]$ at powers of $y$.

\begin{theorem}\label{quartic}
Let $X$ be a quartic surface with intersection matrix $\left(\begin{smallmatrix}-2 & \,\,3 \\\,\,3 & 0 \\ \end{smallmatrix}\right)$. Then the Cox ring of $X$ is isomorphic to the $k-$algebra
\begin{equation}\label{moriCoxRing}
k[x,y_1,y_2,z_1,z_2,t]/(h_1,h_2)
\end{equation}where $\deg z=\f_1, \deg{y_i}=\f_2, \deg{z_i}=\f_1+\f_2, \deg{t}=3\f_1+2\f_2$. The ideal is generated by two relations $h_1,h_2$ of degree $3\f_1+3\f_2$.
\end{theorem}

\begin{proof}
The sections $x,y_1,y_2,z_1,z_2,t$ generate $\Cox(X)$ by Theorem \ref{generators} by the above discussion. Let $H=\f_1+\f_2$ be the hyperplane divisor of $X$ in the embedding of $X$ as a quartic surface. The fact that there are two minimal relations in degree $3H=3\f_1+3\f_2$ comes from the fact that $R(X,H)$ is generated in degree 1, and so it must be possible to write the sections $ty_1$ and $ty_2$ in terms of  $xy_1,xy_2,z_1,z_2$. It follows that we have at least two relations of the form
$$
h_i=ty_i-f_i(xy_1,xy_2,z_1,z_2)=0.
$$To show that these generate all the relations in $\Cox(X)$, it is sufficient to show that $(h_1,h_2)$ is a prime ideal, since then the ring defined by \eqref{moriCoxRing} is an integral domain that surjects onto $\Cox(X)$, and hence is isomorphic to $\Cox(X)$ since it has Krull dimension 4.

To prove this, note that the polynomial $F=xy_1h_2-xy_2h_1$ is a relation of degree $4H$ which is a polynomial in $xy_1,xy_2,z_1,z_2$. It follows that $F$ is the pullback of the quartic polynomial defining $X$ under the embedding given by $H$. In particular, the elimination ideal $(xy_1h_2-xy_2h_1)=(h_1,h_2)\cap k[x,y_i,z_i]$ is prime. Also, it is straightforward to check that $y_1$ is not a zero-divisor modulo $(h_1,h_2)$, so by the trick quoted above, we conclude that also $(h_1,h_2)$ is prime. \end{proof}

\subsection{A quartic surface containing two plane conics}

Let $X$ be a quartic surface such that a hyperplane section of $X$ splits into two plane conics $\f_1,\f_2$. If $X$ is general it has Picard number 2 with $\Pic(X)=\ZZ\f_1\oplus \ZZ\f_2$ with $\f_i^2=-2, \f_1\cdot \f_2=4$.  The effective cone is generated by $\f_1$ and $\f_2$ and these are the only $(-2)$-curves. The nef cone is generated by $2\f_1+\f_2$ and $\f_1+2\f_2$. Moreover, there are no hyperelliptic divisor classes on $X$, by Proposition \ref{nefp}$iv)$. It is also easy to check that for every ample divisor $D=a\f_1+b\f_2$ with $a\ge 3$ or $b\ge 3$, we have $H^1(X,D-2H)=0$ and $H^2(D-3H)=0$. By Theorem \ref{generators}, these observations show that the degrees of the minimal generators of $\Cox(X)$ are contained in $$\G=\{\f_1,\f_2,\f_1+\f_2,\f_1+2\f_2,2\f_1+\f_2\}.$$ The generators are chosen as follows:
$$
\begin{array}{ll}
H^0(X,\f_1)=\langle x \rangle & H^0(X,\f_1+\f_2)=\langle xy,z_1,z_2,z_3\rangle\\
H^0(X,\f_2)=\langle y\rangle &H^0(X,\f_1+2\f_2)=\langle xy^2, yz_1, yz_2, yz_3, v\rangle\\
&H^0(X,2\f_1+\f_2)=\langle  x^2y, xz_1, xz_2, xz_3, w\rangle
\end{array}
$$

\begin{proposition}
Let $X$ be a quartic surface with $\Pic(X)=\ZZ\f_1\oplus \ZZ\f_2$ with $\f_i^2=-2, \f_1\cdot \f_2=4$. Then the Cox ring of $X$ is given by
$$
\Cox(X)=k[x,y,z_1,z_2,z_3,v,w]/I
$$where $I=(xv-f,yw-g,vw-h)$ is a complete intersection.
\end{proposition}

\begin{proof}
The sections $x,y,z_1,z_2,z_3,v,w$ generate the Cox ring by Theorem \ref{generators}. The three relations come from the fact that $\oplus_{m\ge 0}H^0(X,mH)$ is generated in degree 1, hence it must be possible to write the sections $xv,yw,vw$ in terms of $xy,z_1,z_2,z_3$ (which form a basis for $H^0(X,H)$). The fact that these relations generate the whole ideal of relations in $\Cox(X)$ comes from the fact that $I\cap k[x,y,z_1,z_2,z_3]$ contains the polynomial $fg-xyh$, which must be a constant multiple of the equation defining $X$ as a quartic surface. Hence $I\cap k[x,y,z_1,z_2,z_3]=(fg-xyh)$ is prime. As in the proof of Theorem \ref{quartic}, it follows that $I$ is prime and hence contains all the defining relations of $\Cox(X)$.
\end{proof}

\subsection{A double cover of $\mathbb{F}_4$}\label{hirz}
Let $X$ be a K3 surface with $\Pic(X)$ generated by $\f_1,\f_2$ where $\f_1,\f_2$ are smooth curves with self intersection $0$ and $-2$ respectively and $\f_1\cdot \f_2=1$. The nef cone of $X$ is generated by $\f_1$ and $2\f_1+\f_2$ and there is a unique hyperelliptic class on $X$, which is given by $4\f_1+2\f_2$. If $D=a\f_1+b(2\f_1+\f_2)$ represents a nef divisor class with $a,b>0$, we have $H^1(X,D-f_1)=H^1(X,D-2\f_2)=0$ and thus in the notation of Theorem \ref{generators}, we have
 
$$\G=\{\f_1,\f_2,4\f_1+2\f_2,6\f_1+4\f_2\}.$$

As in the previous examples, we see that the Cox ring of $X$ is given by a quotient
$$
\Cox(X)=k[x,y_1,y_2,z,w]/F
$$where $\deg x=\f_2$, $\deg y_i=\f_1$, $\deg z=4\f_1+2\f_2, \deg w=6\f_1+3\f_2$ and $\deg F=12\f_1+6\f_2$. 

The relation $F$ comes from the fact that the section ring $R(X,4\f_1+2\f_1)$ is generated in degrees 1 and 2 and so it must be possible to write $w^2\in H^0(X,12\f_1+6\f_2)$ in terms of other monomials. In fact, we can understand this relation geometrically if we examine the morphism given by $D=4\f_1+2\f_2$. This divisor is hyperelliptic and therefore defines a double cover $\phi_D:X\to S\subset \PP^5$, where $S\simeq \mathbb{F}_4$ is a Hirzebruch surface. If we choose sections so that $x,y_1,y_2,z$ generate the Cox ring of $S$ (which is a polynomial ring since $S$ is a toric variety) and $w$ defines the ramification divisor of this cover, we see that the relation $F$ will be other form $F=w^2-f(x,y_1,y_2,z)$. Cox rings of K3 surfaces occuring in this manner were studied more generally in \cite[Section 4]{AHL09}.

\section{Doubly elliptic K3 surfaces}\label{2ell}

Consider a K3 surface $X$ with $\Pic(X)$ generated by the classes of two smooth elliptic curves $\f_1,\f_2$. The intersection form on $\Pic(X)$ is given by the matrix $
\left(\begin{smallmatrix}0 & d\\
d & 0
\end{smallmatrix}\right)
$. Cox rings of K3 surfaces of this type were also studied in \cite{AHL09}, where the authors find generators for $\Cox(X)$ and show that the defining ideal is generated by quadrics. In this section we extend this result using classical results on linear systems on such surfaces. In particular, we are able to study the higher syzygies of $\Cox(X)$.

We will in following consider the case $d\ge 2$. If $d=1$, then either $\f_1-\f_2$ or $\f_2-\f_1$ represents an effective $(-2)$-curve and $X$ is the K3 surface studied in the previous section.

When $d\ge2$ it is easy to check that there are no $(-2)$-curves on $X$ and that $\f_1$ and $\f_2$ generate the effective cone of $X$. Using Riemann--Roch, we see that the sections of $\f_1,\f_2$ define base-point free pencils $X\to \PP^1$, so $X$ is an elliptic fibration in two different ways. Also, since $\f_1$ and $\f_2$ are nef, this means that every effective divisor $D=a\f_1+b\f_2$ with $a,b\ge 0$ is nef and ample if $a,b>0$.

We fix independent sections $x_1,x_2 \in H^0(X,\f_1)$ and $y_1,y_2\in H^0(X,\f_2)$. These give four linearly independent monomials $x_iy_j$ in degree $H=\f_1+\f_2$. If $d=2$, these monomials span $H^0(X,H)$, and when $d\ge 3$, we add $d-2$ sections $z_1,\ldots,z_{d-2}$ for a basis. In the latter case, the algebra $R(X,H)$ is generated in degree 1, while when $d=2$, $H$ is hyperelliptic and we must add a generator $z$ of degree $2H$, in accordance with Proposition \ref{nefp}v).

\begin{lemma}\label{generated}
Let $d\ge2$ and $X$ be a K3 surface with intersection matrix $
\left(\begin{smallmatrix}0 & d\\
d & 0
\end{smallmatrix}\right)$and let $H=\f_1+\f_2$.
\begin{itemize}
\item If $d=2$, $\Cox(X)$ is generated by $x_1,x_2,y_1,y_2,z$
\item If $d\ge 3$, $\Cox(X)$ is generated by $x_1,x_2,y_1,y_2,z_1,\ldots,z_{d-2}$
\end{itemize} where $\deg x_i=\f_1, \deg y_i=\f_2, \deg z=2H$ and $\deg z_i=H$.
\end{lemma}
\begin{proof}
This is immediate from Theorem \ref{generators} and the discussion above, by noting that $\f_1+\f_2$ is hyperelliptic only if $d=2$, that there are no $(-2)$-curves on $X$, and the fact that when $D=a\f_1+b\f_2$ is a nef divisor with $a>1$ or $b>1$, either $E=\f_1$ or $E=\f_2$ gives $H^1(X,D-2E)=0$.
\end{proof}

\noindent {\bf Example.} When $d=2$, the divisor class $H=\f_1+\f_2$ is hyperelliptic and the sections $x_1y_1,x_1y_2,x_2y_1,x_2y_2$ define a morphism $\phi_H:X\to \mathbb{P}^3$ onto the quadric surface $Q=\{u_1u_4-u_2u_3=0\}$. The branch locus of $\phi$ is a curve of bidegree $(4,4)$ on $Q$. By the lemma above the Cox ring of $X$ is a quotient of the multigraded polynomial ring $k[x_1,x_2,y_1,y_2,z]$. Since the Cox ring has dimension 4, the defining ideal is generated by exactly one relation. To see what this relation is, we note that by Proposition \ref{koszulK3} the multiplication map  $H^0(X,3H)\otimes H^0(X,H)\to H^0(X,4H)$ is surjective, which means  that we can write the section $z^2$ in terms of the other 34 monomials $x_iy_jz$ of degree $4H$. \qed

\medskip

\noindent In the following we will assume $d\ge 3$, in which case the divisor $H$ is very ample (e.g., by \cite{SD74}). We will consider the polynomial ring \begin{equation}\label{R}R=k[x_1,x_2,y_1,y_2,z_1,\ldots,z_{d-2}]\end{equation} with the multigrading $\deg x_i=\f_1, \deg y_i=\f_2, \deg z_i=\f_1+\f_2$. By Lemma \ref{generated}, $R$ surjects onto  $\Cox(X).$ The induced  $\mathbb{G}^2_m$-action on $R$ is given by $$(t_1,t_2)\cdot x_i=t_1x_i,\quad (t_1,t_2)\cdot y_i=t_2y_i, \mbox{ and } (t_1,t_2)\cdot z_i=t_1t_2z_i.$$It is easy to check, using toric geometry \cite{Cox95} ,that the toric variety $Y=\Spec R \gq \mathbb{G}^2_m$ is isomorphic to a rank 4 quadric in $\PP^{d+1}$, which is a rational normal scroll. This gives us an embedding of $X$ into $Y$.

 In fact, we can make this embedding a little more explicit. Consider $X$ as embedded as a surface of degree $2d$ in $\PP^{d+1}$ using sections of degree $H$. Taking the sections $u_{ij}=x_iy_j, u_i=z_i \in H^0(X,H)$, we see that $X$ lies on the rank 4 quadric $$Y=Z(u_{11}u_{22}-u_{12}u_{21}).$$The two rulings of this scroll cut out the linear systems $|\f_1|$ and $|\f_2|$ and where the general fibers are elliptic normal curves of degree $d$.

The point of all this is that the sections $x_i,y_i,z_i$ define Cox coordinates on the quadric $Y$ and relations in the ideal of $I_{X|Y}$ correspond to relations in the Cox ring. We will use this observation to describe the defining relations in $\Cox(X)$ geometrically. For example, if $d=3$, we find that $X$ is embedded as a degree 6 complete intersection of the above rank 4 quadric and a cubic hypersurface $Y$ in $\PP^4$. The Cox ring of $X$ therefore has a single defining relation $F=0$ where $F$ is the bihomogeneous polynomial of degree $3H$ defining $X$ in $Y$.

Consider the minimal multigraded resolution of $\Cox(X)$ as an $R$-module:
\begin{equation}\label{resol}\cdots \to \bigoplus_D  R(-D)^{b_{2,D}} \rightarrow\bigoplus_D  R(-D)^{b_{1,D}}\to R\to R/I \to 0.\end{equation}
Using the embedding of $X$ into $Y$, we will see below that $\Cox(X)$ has Betti numbers equal to those of elliptic normal curves. To prove this, we need one preparatory lemma:

\begin{lemma}\cite[Lemma 4.1]{Sc86}\label{gorenstein}
For a $0$-dimensional non-degenerate subscheme $Z$ of $\PP^{d-2}$ of degree $d$, the following are equivalent:
\begin{enumerate}
\item The homogeneous coordinate ring of $D$ is Gorenstein. 
\item $\O_Z$ has an $\O=\O_{\PP^{d-2}}$-module resolution of type
$$
0\to \O(-d)\to\O(-d+2)^{\beta_{d-3}}\to \ldots\to  \O(-2)^{\beta_{1}} \to \O \to \O_Z \to 0
$$with $\beta_i=i{d-1\choose i+1}-{d-2\choose i-1}$.
\item No subscheme $E\subset Z$ of degree $d-1$ is contained in a hyperplane.
\end{enumerate}
\end{lemma}

Consider the elliptic fibration $X\to \PP^1$ given by $\f_1$ and let ${E_\lambda}$, $\lambda\in \PP^1$ be any fiber. We consider ${E_\lambda}$ as subscheme of the ruling given by the rank 4 quadric, and abusing notation, we write ${E_\lambda}\subset  \PP_\lambda^{d-1}$. The lemma above enters when finding the Betti numbers of $E_\lambda$ as an $\O_{\PP_\lambda^{d-1}}$-module: By \cite[Prop. 8.8]{TJ04}, ${E_\lambda}$ is arithmetically normal, hence Cohen-Macaulay, so the Betti numbers of $\O_{E_\lambda}$ are equal to that of a general hyperplane section. Taking a hyperplane $H\simeq \PP^{d-2}$ general so that $H \cap X$ is a smooth canonical curve $C$ and $H$ does not contain any of the fibres in $Y$, we get a 0-dimensional scheme $D$ in $H\simeq\PP^{d-2}$ of degree $d$ and a complete base-point free pencil $|D|$ on $C$. Now we verify that the Betti numbers of $D$ are precisely of the form stated in the lemma. 

To see this, it suffices to verify $(3)$. Suppose to the contrary that $(3)$ does not hold, so that there is a subscheme $F\subset D$ of degree $d-1$ with linear span of dimension $ \le d-3$. Then by the geometric form of Riemann-Roch, $|F|$ is also of dimension 1. But then the remaining point $P=D-F$ is a base-point of $|D|$, a contradiction.

We now state our main theorem of this section.

\begin{theorem}Let $X$ be a a K3 surface with intersection matrix $\left(\begin{smallmatrix} 0\,\, & d\\ d\,\, & 0\end{smallmatrix}\right)$ for $d\ge 4$. Consider $X$ embedded in the rank 4 quadric $Y$, whose Cox ring $R=Cox(Y)$ is the polynomial ring in \eqref{R}. Then $\Cox(X)$ has a minimal resolution of the form

{\footnotesize
\begin{equation*}
0   \hpil  R_{}(-dH)  \hpil R_{}((-d+2)H )^{\beta _{d-3}}   \hpil    \cdots   \hpil   R_{}(-2H)^{\beta _1}  \hpil 
R_{} \hpil \Cox(X) \hpil 0,
\end{equation*}}where ${\beta}_i = i{d-1\choose i+1}-{d-2\choose i-1}$. In particular, the defining ideal is generated in degree $2\f_1+2\f_2$.\end{theorem}

\begin{proof}
By sheaf-$\Cox(Y)$-module correspondence on a toric variety \cite{Cox95}, we find that \eqref{resol} gives the minimal resolution of $\O_X$:
\begin{equation}\label{resol2}
\cdots \to \bigoplus_D  \O_Y(-D)^{\beta_{2,D}} \rightarrow\bigoplus_D  \O_Y(-D)^{\beta_{1,D}}\to \O_Y\to \O_X \to 0
\end{equation}We use the method of Schreyer \cite{Sc86} to determine the Betti numbers of the above resolution. Let ${E_\lambda}$ be any member of the linear system $|\f_1|\simeq \PP^1$. Restriction of \eqref{resol2}  to the fibre containing ${E_\lambda}$ gives a map of complexes

$$\begin{array}[c]{ccccccccc}
\cdots  \rightarrow &  \bigoplus_D \O_Y(-D)^{\beta_{2,D}}\rightarrow & \bigoplus_D \O_Y(-D)^{\beta_{1,D}}& \rightarrow & \O_Y& \rightarrow & \O_X \rightarrow& 0\\
&  \downarrow &  \downarrow &  &\downarrow &  & \downarrow &\\
\cdots  \rightarrow &   \O_{\PP^{d-1}}(-3)^{\alpha_{2}}\rightarrow & \O_{\PP^{d-1}}(-2)^{\alpha_{1}}& \rightarrow & \O_{\PP^{d-1}}& \rightarrow & \O_{E_\lambda} \rightarrow& 0\\
\end{array}$$
with $\alpha_i= i{d-1\choose i+1}-{d-2\choose i-1}$. By \cite[Thm. 3.1]{Sc86}, the minimal resolution of $\O_{E_\lambda}$ as an $\O_{\PP^{d-1}}$-module lifts to a minimal resolution of $\O_X$ as an $\O_Y$-module, provided that the Betti numbers $\beta_\lambda$ are the same for all $\lambda\in \PP^1$. But this condition is clear from the discussion above. Now, since the syzygies of $\Cox(X)$ restrict to syzygies of $\O_{E_\lambda}$, and since we may repeat the argument with $\f_2$ instead of $\f_1$, we see that the resolution must be of the form above.  The proof is complete.\end{proof}

\medskip

\noindent {\bf Example.} When $d=5$, $H$ embeds $X$ as a degree 10 surface in $\PP^6$. If $X$ is general, it is known from results by Mukai that $X$ is a 3-fold linear section of the Grassmannian $G(2,5)$ and the rank 4 quadric. It follows that the ideal $I_X$ is generated by the five maximal pfaffians of a $5\times 5$ matrix and the quadric $Q$. This gives five relations in the Cox ring and the resolution is given by
$$
0\to R(-5H)\to R(-3H)^5\to R(-2H)^5\to R\to \Cox(X)\to 0.
$$

\bigskip

\noindent\emph{Acknowledgement.} I wish to thank Prof. Kristian Ranestad for his continuous encouragement and advice. I also want to thank Antonio Laface for useful conversations.

\end{document}